\documentclass[10pt,twoside,reqno]{amsart}
\usepackage{mathrsfs}
\usepackage{csquotes}
\usepackage{mathtools}
\usepackage{float}
\usepackage[margin=3cm]{caption}
\usepackage[margin=1.45in]{geometry}
\usepackage{url}

\newcommand\numberthis{\addtocounter{equation}{1}\tag{\theequation}}

\newtheorem{theorem}{Theorem}[]
\newtheorem{corollary}[]{Corollary}
\newtheorem{lemma}[]{Lemma}

\newtheorem*{primepairconjecture}{Prime Pair Error Conjecture}
\newtheorem*{primepairlowerconjecture}{Prime Pair Error Lower Bound  Conjecture}
\newtheorem*{Vaughanconjecture}{Vaughan's Conjecture}
\newtheorem*{Vaughanlower}{Vaughan Lower Bound Conjecture}
\newtheorem*{acknowledgments}{Acknowledgments}

\makeatletter
\newcommand{\tpmod}[1]{{\@displayfalse\pmod{#1}}}
\makeatother

\begin{document}	

\title[The error term in counting prime pairs]{The error term in counting prime pairs}

\author[L. Chou]{Leon Chou}

\address{Department of Mathematics, University of California, Berkeley, 970 Evans Hall, Berkeley, California 94720, United States of America}

\email{knoeleon6@berkeley.edu}

\author[S. Haag]{Summer Haag}

\address{Department of Mathematics, University of Colorado Boulder, 2300 Colorado Avenue, Boulder, Colorado 80309, United States of America}

\email{suha3163@colorado.edu}

\author[J. Huryn]{Jake Huryn}

\address{Department of Mathematics, The Ohio State University, 100 Math Tower, 231 West 18th Avenue, Columbus, Ohio 43210, United States of America}

\email{huryn.5@osu.edu}

\author[A. Ledoan]{Andrew Ledoan}

\address{Department of Mathematical Sciences, Northern Illinois University, 320 Watson Hall,  1425 West Lincoln Highway, DeKalb, Illinois 60115, United States of America}

\email{aledoan@niu.edu}

\makeatletter
\@namedef{subjclassname@2020}{
\textup{2020} Mathematics Subject Classification}
\makeatother



\begin{abstract}
We relate the size of the error term in the Hardy--Littlewood conjectured formula for the number of prime pairs to the $L^{1}$ norm of an exponential sum over the primes formed with the von Mangoldt function.
\end{abstract}

\maketitle

\thispagestyle{empty}

\section{Introduction} \label{section 1}

Let $S (\alpha)$ denote the exponential sum
\[
S (\alpha)
 = \sum_{n \leq N} \Lambda (n) e (n \alpha),
\]
where $e (u) = e^{2 \pi i u}$ and $\Lambda (n)$ is von Mangoldt's function defined by
\[
\Lambda (n)
 = \left\{ \begin{array}{ll}
      \log p & \mbox{if $n = p^{m}$ for some prime $p$ and integer $m \geq 1$,} \\
      0 & \mbox{otherwise.}
\end{array}
\right.
\]
The usual generating function for prime pairs is
\begin{equation} \label{eq-1}
\lvert S (\alpha) \rvert^{2}
 = \sum_{\lvert k \rvert \leq N} \psi_{2} (N, k) e (k \alpha),
\end{equation}
where
\begin{equation} \label{eq-2}
\psi_{2} (N, k)
 = \sum_{\substack{n, n^{\prime} \leq N \\ n^{\prime} - n = k}} \Lambda(n) \Lambda (n^{\prime}).
\end{equation}
If $k$ is zero, then the prime number theorem gives
\begin{equation} \label{eq-3}
\psi_{2} (N, 0)
 = \sum_{n \leq N} \Lambda(n)^{2}
 = N \log N + O (N) \hspace{1em} \textrm{as $N \to \infty$.}
\end{equation}
If $k$ is nonzero and even, then the Hardy--Littlewood conjecture \cite{HardyLittlewood1922} for prime pairs states that
\[
\psi_{2} (N, k)
 \sim \mathfrak{S} (k) (N - \lvert k \rvert) \hspace{1em} \textrm{as $N \to \infty$,}
\]
where
\[
\mathfrak{S} (k)
 = \left\{ \begin{array}{ll}
      {\displaystyle 2 C_{2} \prod_{\substack{ p \mid k \\ p > 2}} \left(\frac{p - 1}{p - 2}\right)} & \mbox{if $k$ is nonzero and even,} \\
      0 & \mbox{if $k$ is  odd,}
\end{array}
\right.
\]
and
\[
C_{2}
 = \prod_{p > 2} \left( 1 - \frac{1}{(p - 1)^{2}}\right)
 = 0.66016\ldots.
\]
If $k$ is odd, then $\psi_{2} (N, k)$ has nonzero terms only when $n$ or $n^{\prime}$ is a power of two. Hence,
\[
\psi_{2} (N, k)
 \ll (\log N)^{2}.
\]

From the recent works of Zhang \cite{Zhang2014}, Maynard \cite{Maynard2013, Maynard2015}, and Tao \cite{Polymath2014, Polymath2015} we now know that there exist infinitely many integers $k$ and a constant $c_{1}$ for which
\[
\psi_{2} (N, k)
 \gg \frac{N}{(\log N)^{c_{1}}},
\]
but the method used is unlikely to lead to asymptotic formulas or to prove this result for all even $k$. In this paper, we are interested with the size of the error in the Hardy--Littlewood conjectured formula for $\psi_{2} (N, k)$. We show, in particular, that it cannot always be small. This problem has previously been dealt with by Montgomery and Vaughan \cite{MontgomeryVaughan1973} for sums of primes. We will see that their method is less successful here. In the first place, for a given fixed even value of $k$, nothing that is nontrivial is currently known. Heath-Brown \cite{Heath-Brown1994} proposed the following problem for twin primes at an Oberwolfach meeting in 1994:

\begin{displayquote}
If $p_{n}$ denotes the $n$th twin prime, then
\[
p_{n + 1} - p_{n}
 = \Omega ((\log p_{n})^{2}).
\]
Improve on this result. Hence, or otherwise, improve on the trivial relation
\[
\pi_{2} (x) - c_{2} \int_{2}^{x} \frac{dt}{(\log t)^{2}}
 = \Omega (1),
\]
where $c_{2}$ is a positive constant.
\end{displayquote}
This is equivalent to proving
\[
\psi_{2} (N, 2)
 \neq \mathfrak{S} (2) N + O ((\log N)^{2}),
\]
but there has been no progress on this problem. Thus, instead of individual $k$, we will consider the total error, or variance, over all $k$ and define
\[
E (N)
 = \sum_{1 \leq \lvert k \rvert \leq N} (\psi_{2} (N, k) - \mathfrak{S} (k) (N - \lvert k \rvert))^{2}.
\]

\begin{figure}[H]
\centering$
\begin{array}{c}
\includegraphics[width=85mm]{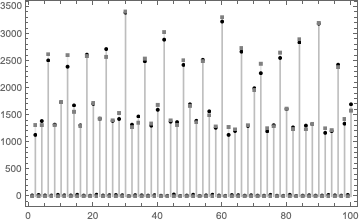}
\end{array}$
\caption{$\psi_{2} (1000, k)$ (black dot) and $\mathfrak{S} (k) (1000 - k)$ (gray square) for $1 \leq k \leq 100$.}
\label{figure 1}
\end{figure}

\begin{figure}[H]
\centering$
\begin{array}{c}
\includegraphics[width=85mm]{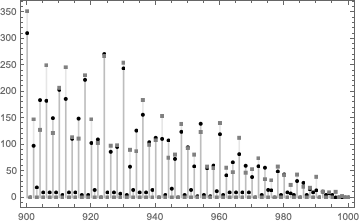}
\end{array}$
\caption{$\psi_{2} (1000, k)$ and $\mathfrak{S} (k) (1000 - k)$ for $900 \leq k \leq 1000$.}
\label{figure 2}
\end{figure}

Figures \ref{figure 1} and \ref{figure 2} give a clear indication of the close agreement between $\psi_{2} (N, k)$ and $\mathfrak{S} (k) (N - \lvert k \rvert)$. Based on this and more extensive numerical evidence supplied by Funkhouser et al. \cite{FunkGoldSengSeng2020}, we propose the following conjecture concerning $E (N)$.

\begin{primepairconjecture}
There exists a positive constant $c_{3}$ such that
\[
E (N) \sim c_{3} N^{2} (\log N)^{2} \hspace{1em} \textrm{as $N \to \infty$.}
\]
\end{primepairconjecture}

We show in Theorem \ref{theorem-1} that the Prime Pair Error Conjecture implies the Riemann Hypothesis. Obtaining upper bounds for $E (N)$ is a very hard problem.  Assuming the Generalized Riemann Hypothesis for Dirichlet's $L$-functions, we have the bound
\[
E (N)
 \ll N^{5 / 2} (\log N)^{c_{4}}
\]
for some positive constant $c_{4}$. This falls short of the bound in the Prime Pair Error Conjecture by a factor of $N^{1 / 2}$.

Our main results connect $E (N)$ to the $L^{1}$ norm of $S (\alpha)$. The latter was first studied by Vaughan \cite{Vaughan1988}. He proved the lower bound
\begin{equation} \label{eq-4}
\int_{0}^{1} \lvert S (\alpha) \rvert \, d \alpha
 \gg N^{1 / 2}.
\end{equation}
It is easy to obtain an upper bound not much larger than this. By the Cauchy--Schwarz inequality and Parseval's equation, together with \eqref{eq-1} through \eqref{eq-3},
\begin{align*}
\int_{0}^{1} \lvert S (\alpha) \rvert \, d \alpha
 & \leq \left(\int_{0}^{1} \lvert S (\alpha) \rvert^{2} \, d \alpha\right)^{1 / 2}
 = \psi_{2} (N, 0)^{1 / 2} \\
 & \leq (1 + o (1)) N^{1 / 2} (\log N)^{1 / 2} \hspace{1em} \textrm{as $N \to \infty$.}
\end{align*}
This upper bound has been improved slightly by Goldston \cite{Goldston2000}, who proved, for any positive $\epsilon$ and $N \geq N_{0} (\epsilon)$,
\begin{equation} \label{eq-5}
\int_{0}^{1} \lvert S (\alpha) \rvert \, d \alpha
 \leq \left[\left(\frac{1}{2} + \epsilon\right) N \log N\right]^{1 / 2}.
\end{equation}
Equations \eqref{eq-4} and \eqref{eq-5} are respectively the current best upper and lower bounds known for the $L^1$ norm of $S (\alpha)$. As we will see below, it appears to be a difficult problem to improve on either of these bounds. ``It seems quite likely,'' Vaughan \cite{Vaughan1988} wrote, that there exists a positive constant $c_{5}$ such that the $L^{1}$ norm of $S (\alpha)$ is asymptotically $c_{5} (N \log N)^{1 / 2}$ as $N \to \infty$, ``but if true this must lie very deep.'' This problem is still unsolved. We will formulate it as the following conjecture.
\begin{Vaughanconjecture}
There exists a positive constant $c_{6}$ satisfying $0 < c_{6} \leq 1 / \sqrt{2}$ such that
\[
\int_{0}^{1} \lvert S (\alpha) \rvert \, d \alpha \sim c_{6} N^{1 / 2} (\log N)^{1 / 2} \hspace{1em} \textrm{as $N \to \infty$.}
\]
\end{Vaughanconjecture}

We emphasize that both the Prime Pair Error Conjecture and Vaughan's Conjecture are extremely difficult problems and there is little hope at present for making progress on proving them. What seems like a much easier problem is the following conjecture.
\begin{primepairlowerconjecture}
We have
\[
E (N)
 \gg N^{2} (\log N)^{2}.
\]
\end{primepairlowerconjecture}

\begin{figure}[H]
\centering$
\begin{array}{c}
\includegraphics[width=85mm]{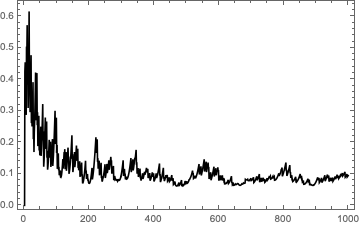}
\end{array}$
\caption{$E (N) / (N^{2} (\log N)^{2})$ for $1 \leq N \leq 1000$.}
\label{figure 3}
\end{figure}

\begin{table}[H]
\begin{center}
\begin{tabular}{|c|c|}
\hline
$N$ & $E (N) / (N^{2} (\log N)^{2})$ \\ 
\hline
$1 \times10^{3}$& $0.09464$ \\
$1 \times 10^{4}$ & $0.12327$ \\
$2 \times 10^{4}$ & $0.13061$ \\
$3 \times 10^{4}$ & $0.14507$ \\
$4 \times 10^{4}$ & $0.15081$ \\
$5 \times 10^{4}$ & $ 0.15480$ \\
$6 \times 10^{4}$ & $0.16124$ \\
$7 \times 10^{4}$ & $0.17745$ \\
$8 \times 10^{4}$ & $0.15953$ \\
$9 \times 10^{4}$ & $0.16192$ \\
$1 \times 10^{5}$ & $0.16857$ \\
\hline
\end{tabular}
\end{center}
\caption{Values of $E (N) / (N^{2} (\log N)^{2})$ truncated to five decimal places.}
\label{table 1}
\end{table}

Figure \ref{figure 3} and Table \ref{table 1} provide some initial numerical evidence in support of the Prime Pair Error Lower Bound Conjecture. We have not been able to prove this conjecture, but there does not appear to be any barrier to doing so. This is somewhat surprising, because Montgomery and Vaughan \cite{MontgomeryVaughan1973} proved essentially this conjecture for the corresponding problem for the error in the Goldbach problem. Moreover, their method which we partly employ here almost works. Many of our results hold when we replace Vaughan's $L^{1}$ Conjecture with the following lower bound conjecture.
\begin{Vaughanlower}
We have
\[
\int_{0}^{1} \lvert S (\alpha) \rvert \, d \alpha \gg N^{1 / 2} (\log N)^{1 / 2} \hspace{1em} \textrm{as $N \to \infty$.}
\]
\end{Vaughanlower}

In Section \ref{section 2}, we will state all of our results in detail. Our first theorem collects together a number of simple results from Corollaries \ref{corollary-1}, \ref{corollary-2}, and \ref{corollary-3}. 

\begin{theorem} \label{theorem-1}
(i) The Prime Pair Error Conjecture implies Vaughan's Lower Bound Conjecture. (ii) If Vaughan's Conjecture holds with $c_{7} < 1 / \sqrt{2}$, then the Prime Pair Error Lower Bound Conjecture is true. (iii) Moreover, if
\[
E (N)
 = o (N^{2} (\log N)^{2} \hspace{1em} \textrm{as $N \to \infty$,}
\]
so that the Prime Pair Error Lower Bound Conjecture is false for all sufficiently large $N$, then Vaughan's Conjecture is true with $c_{7} = 1 / \sqrt{2}$.
\end{theorem}

\section{Statement of results} \label{section 2}

It is not hard to relate $E (N)$ to the classical results on the error in the prime number theorem. This follows Montgomery and Vaughan's work \cite{MontgomeryVaughan1973}, but here it is particularly simple. We easily obtain the following result.

\begin{theorem} \label{theorem-2}
We have
\[
E (N)
 = \Omega (N^{2} (\log \log \log N)^{2}).
\]
If $\Theta$ is the supremum of the real parts of the zeros of Riemann's zeta function, then
\[
E (N)
 = \Omega (N^{1 + 2 \Theta - \epsilon})
\]
for every positive $\epsilon$.
\end{theorem}

If the Riemann Hypothesis is false, then the latter estimate in Theorem \ref{theorem-2} is stronger than the former estimate. Moreover, we see that the Prime Pair Error Conjecture implies that $\Theta = 1 / 2$ and thus the Riemann Hypothesis. Theorem \ref{theorem-2} uses the variability in the number of primes up to $N$ to directly generate the growth in $E (N)$. However, we expect that the growth of $E (N)$ is more strongly influenced by the additive generating function for prime pairs. 

\begin{theorem} \label{theorem-3}
If, for sufficiently large $N$ and a computable large positive constant $c_{8}$,
\begin{equation} \label{eq-6}
\left(\int_{0}^{1} \lvert S (\alpha) \rvert \, d \alpha\right)^{2}
 \leq \frac{1}{2} N \log N - c_{8} N,
\end{equation}
then
\begin{equation} \label{eq-7}
E (N)
 \gg \left[\frac{1}{2} N \log N - \left(\int_{0}^{1} \lvert S (\alpha) \rvert \, d \alpha\right)^2\right]^{2}.
\end{equation}
In particular, \eqref{eq-6} implies
\[
E (N) \gg N^{2}.
\]
\end{theorem}
Unlike the lower bound in Theorem 2, Theorem 3 shows that, subject to \eqref{eq-6}, $E (N)\gg N^{2}$ is never small. We immediately obtain the following corollary.

\begin{corollary} \label{corollary-1}
If, for sufficiently large $N$ and a fixed constant $\delta$ satisfying $0 < \delta < 1 / 2$,
\[
\left(\int_{0}^{1} \lvert S (\alpha) \rvert \, d \alpha\right)^{2}
 \leq \left(\frac{1}{2} - \delta\right) N \log N,
\]
then
\[
E (N)
 \gg \delta^{2} N^{2} (\log N)^{2}.
\]
Hence, if the Vaughan $L^{1}$ Conjecture holds for any constant $c_{5} < 1 / \sqrt{2}$, then the Prime Pair Error Lower Bound Conjecture holds.  
\end{corollary}

\begin{corollary} \label{corollary-2}
If
\[
E (N)
 = o (N^{2} (\log N)^{2}),
\]
then
\[
\int_{0}^{1} \lvert S (\alpha) \rvert \, d \alpha
 \sim \frac{1}{\sqrt{2}} N^{1 / 2} (\log N)^{1 / 2} \hspace{1em} \textrm{as $N \to \infty$.}
\]
\end{corollary}
To see this, we first note that, by Theorem 3, this follows from \eqref{eq-7} subject to \eqref{eq-6}. On the other hand, if \eqref{eq-6} is false for some set of values of $N$, then by Theorem 5 below, for those values of $N$,
\[
\left( \frac{1}{2} N \log N - c_{9} N \right)^{1 / 2}
 \leq \int_{0}^{1} \lvert S (\alpha) \rvert \, d \alpha
 \leq \left(\frac{1}{2} N \log N\right)^{1 / 2} + c_{10} \left(\frac{N}{\log N}\right)^{1 / 2} \log \log N.
\]
The conclusion of Corollary 2 follows independent of the hypothesis that
\[
E (N)
 = o (N^{2}(\log N)^{2}.
\]

Our second main result shows that $E (N)$ and the $L^{1}$ norm of $S (\alpha)$ cannot both be too small.

\begin{theorem} \label{theorem-4}
We have
\[
(E (N) + N^{2} (\log N)^{2}) \left(\int_{0}^{1} \lvert S (\alpha) \rvert \, d \alpha + O (N^{1 / 2} (\log \log N)^{1 / 2})\right)^{2}
 \gg N^{3} (\log N)^{3}.
\]
\end{theorem}

The following corollaries are immediate.

\begin{corollary} \label{corollary-3}
If
\[
E (N)
 \ll N^{2} (\log N)^{c_{11}}
\]
for some positive constant $c_{11}$ satisfying $2 \leq c_{11} < 3$, then
\[
\int_{0}^{1} \lvert S (\alpha) \rvert \, d \alpha
 \gg N^{2} (\log N)^{(3 - c_{11}) / 2}.
\]
In particular, if $c_{11} = 2$ and if
\[
E (N)
 \ll N^{2} (\log N)^{2},
\]
then
\[
\int_{0}^{1} \lvert S (\alpha) \rvert \, d \alpha
 \gg N^{1 / 2} (\log N)^{1 / 2}.
\]
Consequently, the Prime Pair Error Conjecture implies Vaughan's Lower Bound Conjecture.
\end{corollary}

\begin{corollary} \label{corollary-5}
If
\[
\int_{0}^{1} \lvert S (\alpha) \rvert \, d \alpha
 \ll N^{1 / 2} (\log N)^{c_{12}}
\]
for some positive constant $c_{12}$ satisfying $0 < c_{12} \leq 1 / 2$, then
\[
E (N)
 \gg N^{2} (\log N)^{3 - 2 c_{12}}.
\]
\end{corollary}

\begin{corollary} \label{corollary-6}
If
\[
\int_{0}^{1} \lvert S (\alpha) \rvert \, d \alpha
 \ll N^{1 / 2} (\log \log N)^{1 / 2},
\]
then
\[
E (N)
 \gg \frac{N^{2} (\log N)^{3}}{\log \log N}.
\]
\end{corollary}

Many of our results depend on improving the upper bound in \eqref{eq-5}. Turning to this question, we prove the following result.
\begin{theorem} \label{L1Sub} We have
\[
\left( \int_{0}^{1} \lvert S (\alpha) \rvert \, d \alpha \right)^{2}
 \leq \frac{1}{2} N \log N + \frac{3}{2} N \log \log N + O (N).
\]
Alternatively, for some positive constant $c_{10}$ and $N$ sufficiently large, 
\[
\int_{0}^{1} \lvert S (\alpha) \rvert \, d \alpha \leq \left(\frac{1}{2} N \log N \right)^{1 / 2} + c_{10} \left(\frac{N}{\log N}\right)^{1 / 2} \log \log N.
\]
\end{theorem}
From Theorem 3, even small improvements in Theorem 5 would have important consequences for $E (N)$.

\section{Proof of Theorem \ref{theorem-2}} \label{section 3}

For positive $N$, using
\[
\psi (N)
 = \sum_{n \leq N} \Lambda (n),
\]
we obtain
\[
\psi (N)^{2}
 = \sum_{n, n^{\prime} \leq N} \Lambda (n) \Lambda (n^{\prime})
 = \sum_{\lvert k \rvert \leq N} \sum_{\substack{n, n^{\prime} \leq N \\ n^{\prime} - n = k}} \Lambda (n) \Lambda (n^{\prime})
 = \sum_{\lvert k \rvert \leq N} \psi_{2} (N, k).
\]
Let
\[
\psi (N)
 = N + R (N).
\]
The prime number theorem is equivalent to
\[
R (N) = o (N) \hspace{1em} \textrm{as $N \to \infty$.}
\]
Since
\[
\psi (N)^{2}
 = N^{2} + 2 N R (N) + R (N)^{2}
\]
and
\[
\sum_{1 \leq \lvert k \rvert \leq N} (N - \lvert k \rvert) \mathfrak{S} (k)
 = N^{2} + O (N \log N) 
\]
(see \cite{Goldston1984, MontgomerySoundararajan2002}), we deduce that
\[
\sum_{1 \leq \lvert k \rvert \leq N} (\psi_{2} (N, k) - (N - \lvert k \rvert) \mathfrak{S} (k))
 = R (N) (2 N + R (N)) + O (N \log N).
\]
By the Cauchy--Schwarz inequality,
\[
\left\lvert \sum_{1 \leq \lvert k \rvert \leq N} (\psi_{2} (N, k) - (N - \lvert k \rvert) \mathfrak{S} (k)) \right\rvert^{2}
 \leq 2 N E (N).
\]
We conclude that
\[
E (N)
 \geq 2 N \left[R (N)\left(1 + \frac{R (N)}{2 N}\right) + O (\log N)\right]^{2}.
\]
By Littlewood's theorem,
\[
R (N)
 = \Omega_{\pm} (N^{1 / 2} \log \log \log N)
\]
(see \cite{MontgomeryVaughan2007}), and the first assertion of Theorem  \ref{theorem-2} follows. If the Riemann Hypothesis is false, we use
\[
R (N)
 = \Omega_{\pm} (N^{\Theta - \epsilon})
\]
for every positive $\epsilon$ (see \cite{Ingham1932, MontgomeryVaughan2007}). In fact, Littlewood proved that
\[
\limsup_{N \to \infty} \frac{\pm R (N)}{N^{1 / 2} \log \log \log N}
 \geq \frac{1}{2},
\]
and so it follows that
\[
\limsup_{N \to \infty} \frac{E (N)}{N^{2} (\log \log \log N)^{2}}
 \geq \frac{1}{2}.
\]

\section{Proof of Theorem \ref{theorem-3}} \label{section 4}

Our main results depend on approximating $\lvert S (\alpha) \rvert^{2}$ by
\[
V_{y} (\alpha)
 = \sum_{q \leq y} \frac{\mu(q)^{2}}{\phi(q)^{2}} \sum_{\substack{1 \leq a \leq q \\ (a, q) = 1}} \left\lvert I \left(\alpha - \frac{a}{q}\right) \right\rvert^{2},
\]
where
\[
I (\alpha)
 = \sum_{n \leq N} e (n \alpha),
\]
where $\mu (q)$ is the M\"{o}bius function and $\phi (q)$ is the Euler totient function. Here, $V_{y} (\alpha)$ is obtained by taking the diagonal terms from squaring the absolute value of the major arc approximation $J_{Q} (\alpha)$ of $S (\alpha)$ used by Goldston \cite{Goldston2000} and removing the restricted support for each arc. (This is explained in Section \ref{section 7}. Our use of $V_{y} (\alpha)$ is independent of these heuristics.) Obviously, $V_{y} (\alpha)$ is real-valued and nonnegative. Using the relation
\[
\lvert I (\alpha) \rvert^{2}
 = \sum_{\lvert k \rvert \leq N} (N - \lvert k \rvert) e (k \alpha),
\]
we obtain the Fourier series
\begin{align*}
V_{y} (\alpha)
 & = \sum_{\lvert k \rvert \leq N} (N - \lvert k \rvert) \left(\sum_{q \leq y} \frac{\mu (q)^{2}}{\phi (q)^{2}} \sum_{\substack{1 \leq a \leq q \\ (a, q) = 1}} e \left(-\frac{k a}{q}\right)\right) e (k \alpha) \\
 & = \sum_{\lvert k \rvert \leq N} (N - \lvert k \rvert) \mathfrak{S}_{y} (k) e (k \alpha),
\end{align*}
where
\[
\mathfrak{S}_{y} (k)
 = \sum_{q \leq y} \frac{\mu(q)^{2}}{\phi(q)^{2}} c_{q} (-k)
\]
and
\[
c_{q} (n)
 = \sum_{\substack{1 \leq a \leq q \\ (a, q) = 1}} e \left(\frac{a n}{q}\right).
\]
The Ramanujan sum $c_{q} (n)$ is an even function in $n$ and a multiplicative function of $q$ for a fixed value of $n$, and satisfies
\begin{equation} \label{eq-8}
c_{q} (n)
 = \sum_{d \mid (q, n)} d \mu \left(\frac{q}{d}\right)
 = \mu \left(\frac{q}{(q, n)}\right) \frac{\phi (q)}{\displaystyle \phi \left(\frac{q}{(q, n)}\right)}
\end{equation}
(see \cite[Theorem 4.1, (4.7)]{MontgomeryVaughan2007}).

It is now easy to verify that $\mathfrak{S} (k)$ in the Hardy--Littlewood conjecture for prime pairs has the series representation
\[
\mathfrak{S} (k)
 = \sum_{q = 1}^{\infty} \frac{\mu (q)^{2}}{\phi (q)^{2}} c_{q} (-k)
\]
for $k \neq 0$. While we do not define $\mathfrak{S} (k)$ for $k = 0$, we have from \cite[2.1.1, Ex. 17]{MontgomeryVaughan2007} that
\begin{equation} \label{eq-9}
\mathfrak{S}_{y} (0)
 = \sum_{q \leq y} \frac{\mu(q)^{2}}{\phi(q)}
 = \log y + O (1).
\end{equation}
Let 
\begin{equation} \label{eq-10}
a_{0} (N, y)
 = \sum_{n \leq N} \Lambda(n)^{2} - N \mathfrak{S}_{y} (0).
\end{equation}
By \eqref{eq-3} and \eqref{eq-9}, for $1 \leq y \leq N$,
\begin{equation} \label{eq-11}
a_{0} (N, y)
 = N \log \left(\frac{N}{y}\right) + O (N).
\end{equation}
Let
\[
V (\alpha)
 = \sum_{1 \leq \lvert k \rvert \leq N} (N - \lvert k \rvert) \mathfrak{S} (k) e (k \alpha).
\]
By Parseval's identity,
\begin{align*}
E (N)
 & = \int_{0}^{1} \left\lvert\lvert S (\alpha) \rvert^{2} - V (\alpha) - \sum_{n \leq N} \Lambda(n)^{2} \right\rvert^{2} d \alpha \\
 & = \int_{0}^{1} \lvert (\lvert S (\alpha) \rvert^{2} - V_{y} (\alpha) - a_{0} (N, y)) - (V (\alpha) - V_{y} (\alpha) + N \mathfrak{S}_{y} (0)) \rvert^{2} \, d \alpha.
\end{align*}

Let now
\[
J
 = J (N, y)
 = \int_{0}^{1} \lvert \lvert S (\alpha) \rvert^{2} - V_{y} (\alpha) - a_{0} (N, y) \rvert^{2} \, d \alpha
\]
and
\[
W
 = W (N, y)
 = \int_{0}^{1} \lvert V (\alpha) - V_{y} (\alpha) + N \mathfrak{S}_{y} (0) \rvert^{2} \, d \alpha .
\]
By the Cauchy--Schwarz inequality,
\begin{equation} \label{eq-12}
E (N)
  \geq J - 2 (J W)^{1 / 2} + W
 = (J^{1 / 2} - W^{1 / 2})^{2}.
\end{equation}
We obtain a lower bound for $E (N)$ from a lower bound for $J$ and an asymptotic formula for $W$. From Parseval's identity and \cite[Theorem 2]{GoldstonHuntsNgotiaoco2017}, for $1 \leq y \leq N^{1 / 2}$,
\begin{equation}
\begin{split} \label{eq-13}
W
 & = \sum_{1 \leq \lvert k \rvert \leq N} (N - \lvert k \rvert)^{2} (\mathfrak{S} (k) - \mathfrak{S}_{y} (k))^{2} \\
 & = \frac{1}{3} N^{3} H (y) - \frac{1}{4} N^{2} \left[\log \left(\frac{N}{y^{2}}\right)^{2}\right] + M N^{2} \log \left(\frac{N}{y^{2}}\right) \\ & \hspace{1em} + O (N^{2}) + O \left(\frac{N^{2}}{y^{1 / 2}} \log (2 N)\right),
\end{split}
\end{equation}
where
\[
H (y)
 = \frac{L}{y^{2}} (1 + o (1)), \hspace{1em}
L
 = \prod_{p} \left(1 + \frac{2 p - 1}{p (p - 1)^{2}}\right),
\]
and
\[
M
 = \frac{3}{4} - \frac{1}{2} \log (2 \pi) + \frac{1}{2} \sum_{p} \frac{(p - 2) \log p}{p (p - 1)^{2}}.
\]
Hence, for  $1 \leq y \leq N^{1 / 2}$,
\begin{equation} \label{eq-14}
W
 \asymp \frac{ N^{3}}{y^{2}}.
\end{equation}
The implied upper and lower bound constants are computable since, as mentioned in \cite{GoldstonHuntsNgotiaoco2017}, the term $O (N^{2})$ in \eqref{eq-13} can be replaced by $c_{13} N^{2} + o (N^{2})$ for a computable constant $c_{13}$.

Montgomery and Vaughan \cite{MontgomeryVaughan1973} obtained a lower bound for what corresponds to $J$ for the Goldbach problem. Their argument is that, by the Cauchy--Schwarz inequality,
\begin{equation} \label{eq-15}
J
 \geq T^{2},
\end{equation}
where
\[
T
 = T (N, y)
 = \int_{0}^{1} \lvert \lvert S (\alpha) \rvert^{2} - V_{y} (\alpha) - a_{0} (N, y) \rvert \, d \alpha.
\]
But, by \eqref{eq-10},
\begin{align*}
T
 & \geq \int_{0}^{1} (\lvert S (\alpha) \rvert^{2} - V_{y} (\alpha) - a_{0} (N, y)) \, d \alpha \\
 & = \int_{0}^{1} \lvert S (\alpha) \rvert^{2} \, d\alpha -  \int_{0}^{1}V_{y} (\alpha) \, d \alpha - a_{0} (N, y) \\
 & = \sum_{n \leq N} \Lambda(n)^{2} - N \mathfrak{S}_{y} (0) - a_{0} (N, y)
 \equiv 0,
\end{align*}
and this argument fails. In the Goldbach problem, the term $a_{0} (N, y) = 0$, which is why it succeeds there.  
We now modify this argument to obtain a lower bound for $J$ in terms of the $L^1$ norm of $S (\alpha)$.

\begin{lemma} \label{lemma-1}
If, for $y_{0} \leq y \leq N^{1 / 2}$ for some $y_{0}$ satisfying $N^{1 / 2} / \log N \leq y_{0} \leq N^{1 / 2}$ and a computable positive constant $c_{14}$,
\[
\left(\int_{0}^{1} \lvert S (\alpha) \rvert \, d \alpha\right)^{2}
 \leq a_{0} (N, y) - c_{14}\frac{N^{3 / 2}}{y},
\]
then
\[
E (N)
 \geq \frac{1}{4} \left[a_{0} (N, y) - \left(\int_{0}^{1} \lvert S (\alpha) \rvert \, d \alpha\right)^{2}\right]^{2}
 \gg \frac{N^{3}}{y^{2}}.
\]
\end{lemma}

Theorem 3 is the special case of Lemma \ref{lemma-1} where we take $y_{0} = N^{1 / 2}$.

\begin{proof}
We utilize the inequality, for any two real numbers $a$ and $b$,
\[
\lvert a^{2} - b^{2} \rvert \geq (\lvert a \rvert - \lvert b \rvert)^{2}.
\]
This inequality holds true because of the identity
\[
\lvert a^{2} - b^{2} \rvert
 = (\lvert a \rvert - \lvert b \rvert)^{2} + 2 \min (\lvert a \rvert, \lvert b \rvert) \lvert\lvert a \rvert - \lvert b \rvert\rvert,
\]
which in turn holds true since we can assume without loss of generality that $a \geq b \geq 0$, in which case it is obvious. Since $V_{y} (\alpha) \geq 0$ and  $a_{0} (N, y) > 0$ for $1 \leq y \leq N^{1 / 2}$, we have
\begin{align*}
T
 & \geq \int_{0}^{1} \lvert \lvert S (\alpha) \rvert^{2} - a_{0} (N, y) \rvert \, d \alpha - \int_{0}^{1} V_{y} (\alpha) \, d \alpha \\
 & \geq \int_{0}^{1} \lvert \lvert S (\alpha) \rvert - a_{0} (N, y)^{\frac{1}{2}} \rvert^{2} \, d \alpha - N \mathfrak{S}_{y} (0) \\
 & = \int_{0}^{1} \lvert S (\alpha) \rvert^{2} \, d \alpha - 2 a_{0} (N, y)^{1 / 2} \int_{0}^{1} \lvert S (\alpha) \rvert \, d \alpha + a_{0} (N, y) - N \mathfrak{S}_{y} (0) \\
 & = 2 a_{0} (N, y)^{1 / 2} \left( a_{0} (N, y)^{1 / 2} - \int_{0}^{1} \lvert S (\alpha) \rvert \, d \alpha \right).
\end{align*}
Provided
\[
\int_{0}^{1} \lvert S (\alpha) \rvert \, d \alpha
 \leq a_{0} (N, y)^{1 / 2},
\]
substitution into \eqref{eq-15} yields
\begin{align*}
J
 &\geq T^{2}
 \geq 4 a_{0} (N, y) \left(a_{0} (N, y)^{1 / 2} - \int_{0}^{1} \lvert S (\alpha) \rvert \, d \alpha\right)^{2} \\
 & = 4 a_{0} (N, y) \frac{\displaystyle \left[a_{0} (N, y) - \left(\int_{0}^{1} \lvert S (\alpha) \rvert \, d \alpha\right)^{2} \right]^{2}}{\displaystyle \left(a_{0} (N, y)^{1 / 2} + \int_{0}^{1} \lvert S (\alpha) \rvert \, d \alpha\right)^{2}} \\
 & \geq 4 a_{0} (N, y)  \frac{\displaystyle \left[a_{0} (N, y) -\left(\int_{0}^{1} \lvert S (\alpha) \rvert \, d \alpha\right)^{2} \right]^{2}}{ [2 a_{0} (N, y)^{1 / 2}]^{2}} \\
 & = \left[a_{0} (N, y) -\left(\int_{0}^{1} \lvert S (\alpha) \rvert \, d \alpha\right)^{2} \right]^{2}.
\end{align*}
Observe that the lower bound for $J$ we have just obtained is $\leq a_{0} (N, y)^{2} \ll N^{2} (\log N)^{2}$. In order for \eqref{eq-12} to give a lower bound for $E (N)$, we need $J$ to be larger than $W$ which, by \eqref{eq-14}, requires $N^{2} (\log N)^{2} \gg N^{3} / y^{2}$ and, therefore, $ y \gg N^{1 / 2} / \log N$.

Assuming that $J \geq 4 W$, we have by \eqref{eq-12} that
\[
E (N) \geq (J^{1 / 2} - W^{1 / 2})^{2}
 \geq \frac{1}{4} J,
\]
and thus we conclude
\[
E (N)
 \geq \frac{1}{4} \left[a_{0} (N, y) - \left(\int_{0}^{1} \lvert S (\alpha) \rvert \, d \alpha\right)^2 \right]^{2},
\] 
provided that
\[
\left(\int_{0}^{1} \lvert S (\alpha) \rvert \, d \alpha\right)^{2} \leq a_{0} (N, y).
\]

By \eqref{eq-13} and the comment following it, we can find a constant $c_{15}$ such that $W \leq c_{15} N^{3} / y^{2}$, and thus the condition $J \geq 4 W$ is satisfied when
\[
\left[a_{0} (N, y) - \left(\int_{0}^{1} \lvert S (\alpha) \rvert \, d \alpha\right)^2 \right]^{2}
 \geq 4 c_{15} \frac{N^{3}}{y^{2}},
\]
which is equivalent to 
\[
\left(\int_{0}^{1} \lvert S (\alpha) \rvert \, d \alpha\right)^{2}
 \leq a_{0} (N, y) - 2 c_{15}^{1 / 2} \frac{N^{3}}{y^{2}}.
\]
This proves Lemma \ref{lemma-1} with $c_{14} = 2 c_{15}^{1 / 2}$.
\end{proof}

\section{Proof of Theorem \ref{theorem-4}} \label{section 5}

If $1 < p < \infty$ and $1 / p + 1 / q = 1$, and if $f \in L^{p}$ and $g \in L^{q}$, then $f g \in L^{1}$ and H\"{o}lder's inequality is that
\[
\int \lvert f g \rvert
 = \| f g \|_{1} \leq \| f \|_{p} \, \| g \|_{q}.
\]
We apply this inequality with $f = \lvert h \rvert^{1 / 3}$ and $g = \lvert h \rvert^{2 / 3}$, and thus
\[
\| h \|_{1}
 \leq \| \lvert h \rvert^{1 / 3} \|_{p} \, \| \lvert h \rvert^{2 / 3} \|_{q}.
\]
If $p = 3 / 2$, then $q = 3$, and we have
\[
\left(\int_{0}^{1} \lvert h (\alpha) \rvert \, d \alpha\right)^{3}
 \leq \left(\int_{0}^{1} \lvert h (\alpha) \rvert^{1 / 2} \, d \alpha\right)^{2} \int_{0}^{1} \lvert h (\alpha) \rvert^{2} \, d \alpha.
\]

We take
\[
h (\alpha)
 = \lvert S (\alpha) \rvert^{2} - V_{y} (\alpha)
\]
and have
\[
\left(\int_{0}^{1} \lvert \lvert S (\alpha) \rvert^{2} - V_{y} (\alpha) \rvert \, d \alpha\right)^{3}
 \leq \left(\int_{0}^{1} \lvert \lvert S (\alpha) \rvert^{2} - V_{y} (\alpha) \rvert^{1 / 2} \, d \alpha \right)^{2} \int_{0}^{1} \lvert \lvert S (\alpha) \rvert^{2} - V_{y} (\alpha) \rvert^{2} \, d \alpha.
\]
We now set $y = N^{1 / 2}$.  Using $V_{y} (\alpha) \geq 0$, and equations \eqref{eq-10} and \eqref{eq-11}, we have
\begin{align*}
\int_{0}^{1} \lvert \lvert S (\alpha) \rvert^{2} - V_{y} (\alpha) \rvert \, d \alpha
 & \geq \int_{0}^{1} \lvert S (\alpha) \rvert^{2} \, d \alpha - \int_{0}^{1} V_{y} (\alpha) \, d \alpha \\
 & = a_{0} (N, y)
 \sim \frac{1}{2} N \log N \hspace{1em} \textrm{as $N \to \infty$,}
\end{align*}
and therefore
\[
\left(\int_{0}^{1} \lvert \lvert S (\alpha) \rvert^{2} - V_{y} (\alpha) \rvert \, d \alpha\right)^{3}
 \gg N^{3} (\log N)^{3}.
\]
By Parseval's identity, and making use of $W\ll N^2$ from \eqref{eq-14}, we have
\begin{align*}
\int_{0}^{1} \lvert \lvert S (\alpha) \rvert^{2} - V_{y} (\alpha) \rvert^{2} \, d \alpha
 & = \sum_{1 \leq \lvert k \rvert \leq N} (\psi_{2} (N, k) - \mathfrak{S}_{y} (k) (N - \lvert k \rvert))^{2} + a_{0} (N, y)^{2} \\
 & \ll E (N) + W + \frac{1}{4} N^{2} (\log N)^{2} (1 + o (1)) \\
 & \leq E (N) + N^{2} (\log N)^{2} .
\end{align*}
The proof will be completed once we prove
\[
\int_{0}^{1} \lvert \lvert S (\alpha) \rvert^{2} - V_{y} (\alpha) \rvert^{1 / 2} \, d \alpha
 \ll \int_{0}^{1} \lvert S (\alpha) \rvert \, d \alpha + O (N^{1 / 2} (\log \log N)^{1 / 2}).
\]
Using the inequality
\[
\lvert \lvert a \rvert \pm \lvert b \rvert \rvert^{1 / 2}
 \leq \lvert a \rvert^{1 / 2} + \lvert b \rvert^{1 / 2},
\]
we obtain
\[
\int_{0}^{1} \lvert \lvert S (\alpha) \rvert^{2} - V_{y} (\alpha) \rvert^{1 / 2} \, d \alpha
 \leq \int_{0}^{1} \lvert S (\alpha) \rvert \, d \alpha + \int_{0}^{1} \lvert V_{y} (\alpha) \rvert^{1 / 2} \, d \alpha,
 \]
with $y = N^{1 / 2}$. It remains only to prove that
\[
\int_{0}^{1} \lvert V_{y} (\alpha) \rvert^{1 / 2} \, d \alpha \ll N^{1 / 2} (\log \log N)^{1 / 2}.
\]

Let now
\[
w (q)
 = \sum_{\substack{1 \leq a \leq q \\ (a, q) = 1}} \frac{\mu (q)^{2}}{\phi (q)^{2}} \left \lvert I \left(\alpha - \frac{a}{q}\right) \right\rvert^{2}.
\]
We have
\[
\int_{0}^{1} \lvert V_{y} (\alpha) \rvert^{1 / 2} \,d \alpha
 = \int_{0}^{1} \left \lvert \sum_{1 \leq q \leq y_{0}} w (q) + \sum_{y_{0} < q \leq y} w (q) \right \rvert^{1 / 2} d \alpha
 \leq I_{1} + I_{2}, \numberthis \label{eq-16}
\]
where
\[
I_{1}
 = \int_{0}^{1} \left \lvert \sum_{1 \leq q \leq y_{0}} w (q) \right \rvert^{1 / 2} d \alpha
\]
and
\[
I_{2}
 = \int_{0}^{1} \left \lvert \sum_{y_{0} < q \leq y} w (q) \right \rvert^{1 / 2} d \alpha.
\]
The integral $I_{1}$ is treated using the inequality
\[
\left(\sum_{n = 1}^{\infty} \lvert a_{n} \rvert\right)^{1 / 2}
 \leq \sum_{n = 1}^{\infty} \lvert a_{n} \rvert^{1 / 2},
\]
which holds true since, on squaring, the left-hand side is just the diagonal terms of the right-hand side. We obtain
\begin{align*}
I_{1}
 & \leq \int_{0}^{1} \sum_{1 \leq q \leq y_{0}} \sum_{\substack{1 \leq a \leq q \\ (a, q) = 1}} \left\lvert \frac{\mu (q)^{2}}{\phi (q)^{2}} \left\lvert I \left(\alpha - \frac{a}{q}\right) \right\rvert^{2} \right\rvert^{1 / 2} d \alpha \\
 & = \int_{0}^{1} \sum_{1 \leq q \leq y_{0}} \sum_{\substack{1 \leq a \leq q \\ (a, q) = 1}} \frac{\mu (q)^{2}}{\phi (q)} \left\lvert I \left(\alpha - \frac{a}{q}\right) \right\rvert d \alpha \\
 & = \sum_{1 \leq q \leq y_{0}} \mu (q)^{2} \int_{0}^{1} \lvert I (\alpha) \rvert \, d \alpha
 \ll y_{0} \log N,
\end{align*}
by the $L^{1}$ estimate for the Dirichlet kernel. Thus, if $y_{0} = N^{1 / 2} / \log N$,
\begin{equation} \label{eq-17}
I_{1}
 \ll N^{1 / 2}.
\end{equation}
The integral $I_{2}$ is treated using the Cauchy--Schwarz inequality and Parseval's identity. We have, using \eqref{eq-9} in the last step,
\begin{align*}
I_{2}
 & \leq \left(\int_{0}^{1} \sum_{y_{0} < q \leq y} \sum_{\substack{1 \leq a \leq q \\ (a, q) = 1}} \frac{\mu (q)^{2}}{\phi (q)^{2}} \left \lvert I \left( \alpha - \frac{a}{q} \right) \right \rvert^{2} d \alpha\right)^{1 / 2} \\
 & = \left(\sum_{y_{0} < q \leq y} \sum_{\substack{1 \leq a \leq q \\ (a, q) = 1}} \frac{\mu (q)^{2}}{\phi (q)^{2}} \int_{0}^{1} \left\lvert I \left(\alpha - \frac{a}{q}\right) \right\rvert^{2} d \alpha \right)^{1 / 2} \\
 & = \left(N \sum_{y_{0} < q \leq y} \frac{\mu (q)^{2}}{\phi (q)}\right)^{1 / 2}
 = \left[N \left( \log \left(\frac{y}{y_{0}}\right) + O (1) \right)\right]^{1 / 2} \\
 & \ll N^{1 / 2} (\log \log N)^{1 / 2}. \numberthis \label{eq-18}
\end{align*}
Substitution of \eqref{eq-17} and \eqref{eq-18} in \eqref{eq-16} yields Theorem \ref{theorem-4}.

\section{First Proof of Theorem 5} \label{section 6}

If
\[
n
 = \prod_{k = 1}^{r} p_{k}^{a_{k}},
\]
then
\[
\log n
 = \sum_{k = 1}^{r} a_{k} \log p_{k}
 = \sum_{d \mid n} \Lambda (d),
\]
and so by M\"{o}bius inversion
\[
\Lambda (n)
 = \sum_{d \mid n} \mu (d) \log (n / d).
\]
Following \cite{Goldston1992}, we approximate $\Lambda(n)$ by
\[
\Lambda_{R} (n)
 = \sum_{\substack{d \mid n \\ d \leq R}} \mu (d) \log \left(\frac{R}{d}\right),
\]
which truncates the number of divisors $d$ of $n$ and smooths the sum. If $1 < n \leq R$, then the condition $d \leq R$ can be dropped, so that $\Lambda_{R} (n) = \Lambda (n)$. Therefore, we can interpret $\Lambda_{R} (n)$ as an approximation of $\Lambda (n)$ which gets better as $R$ increases up to $n$ where it becomes $\Lambda (n)$.

We thus introduce the approximation for $S (\alpha)$,
\[
S_{R} (\alpha)
 = \sum_{n \leq N} \Lambda_{R} (n) e (n \alpha).
\]
We make use of the formula, for $1\leq R\leq N$,
\[
\sum_{n \leq N} \Lambda_{R} (n)^{2}
 = N \log R + O (N)
\]
due to Graham \cite{Graham1978}. In the range $1\leq R \leq N^{1 / 2} / \log N$, this can be obtained by an elementary calculation. By the prime number theorem, we also obtain, for $1\leq R \leq N$ and any positive constant $c_{13}$,
\[
\sum_{n \leq N} \Lambda_R(n)\Lambda(n)
 = N \log R + O (R)+ O (N (\log N)^{-c_{13}})
\]
(see \cite{Goldston1992}).
By Parseval's identity, we have by \eqref{eq-3} and the above estimates, for $1 \leq R \leq N$,
\begin{align*}
\int_{0}^{1} \lvert S (\alpha) - S_{R} (\alpha) \rvert^{2} \, d \alpha
 & = \int_{0}^{1} \left\lvert \sum_{n \leq N} (\Lambda (n) - \Lambda_{R} (n)) e (n \alpha) \right\rvert^{2} \, d \alpha \\
 & = \sum_{n \leq N} (\Lambda (n) - \Lambda_{R} (n))^{2} \\
 & = \sum_{n \leq N} \Lambda (n)^{2} - 2 \sum_{n \leq N} \Lambda (n) \Lambda_{R} (n) + \sum_{n \leq N} \Lambda_{R} (n)^{2} \\
 & = N \log N - 2 N \log R + N \log R + O(N) \\
 & = N \log \left(\frac{N}{R}\right) + O (N).
\end{align*}
Hence, by the Cauchy--Schwarz inequality, for $1 \leq R \leq N$,
\begin{align*}
\int_{0}^{1} \lvert S (\alpha) - S_{R} (\alpha) \rvert \, d \alpha
 & \leq \left(\int_{0}^{1} \lvert S (\alpha) - S_{R} (\alpha) \rvert^{2} \, d \alpha\right)^{1 / 2} \\
 & \leq \left(N \log \left(\frac{N}{R}\right) + O (N)\right)^{1 / 2}.
\end{align*}
By the reverse triangle inequality,
\begin{align*}
\left \lvert \int_{0}^{1} \lvert S (\alpha) \rvert \, d \alpha - \int_{0}^{1} \lvert S_{R} (\alpha) \rvert \, d \alpha \right \rvert
 & \leq \int_{0}^{1} \lvert S (\alpha) - S_{R} (\alpha) \rvert \, d \alpha \\
 & \leq \left(N \log \left(\frac{N}{R}\right) + O (N)\right)^{1 / 2}.
\end{align*}

We now make use of the following result that we will prove below.

\begin{lemma} \label{lemma-2}
We have
\[
\int_{0}^{1} \lvert S_{R} (\alpha) \rvert \, d \alpha
 \ll R \log N.
\]
\end{lemma}
Thus, we conclude
\[
\int_{0}^{1} \lvert S (\alpha) \rvert \, d \alpha
 \leq \left(N \log \left(\frac{N}{R}\right) + O (N)\right)^{1 / 2} + O (R \log N).
\]
We take  $R = N^{1 / 2} / (\log N)^{3 / 2}$ and on squaring obtain
\[
\left( \int_{0}^{1} \lvert S (\alpha) \rvert \, d \alpha \right)^2
 \leq \frac{1}{2} N \log N + \frac{3}{2} N \log \log N + O (N),
\]
from which
\[
\int_{0}^{1} \lvert S (\alpha) \rvert \, d \alpha
 \leq \left(\frac{1}{2} N \log N \right)^{1 / 2} + c_{10} \left(\frac{N}{\log N}\right)^{1 / 2} \log \log N
\]
for a positive constant $c_{10}$. This proves Theorem 5.

It remains to prove Lemma \ref{lemma-2}. We use
\[
\sum_{n \leq N} e (n \alpha)
 \ll \min \left(N, \frac{1}{\| \alpha \|}\right)
\]
to obtain
\begin{align*}
S_{R} (\alpha)
 & = \sum_{n \leq N} \sum_{\substack{d \mid n \\ d \leq R}} \mu (d) \log \left(\frac{R}{d}\right) e (n \alpha) \\
 & = \sum_{d \leq R} \mu (d) \log \left(\frac{R}{d}\right) \sum_{m \leq N / d} e (d m \alpha) \\
 & \ll \sum_{d \leq R} \log \left(\frac{R}{d}\right) \min \left(\frac{N}{d}, \frac{1}{\| d \alpha \|}\right).
\end{align*}
Thus,
\[
\int_{0}^{1} \lvert S_{R} (\alpha) \rvert \, d \alpha
 \ll \sum_{d \leq R} \log \left(\frac{R}{d}\right) \int_{0}^{1} \min \left(\frac{N}{d}, \frac{1}{\| d \alpha \|}\right) \, d \alpha.
\]
As a function of $\alpha$, the integrand on the right-hand side is a function with period $1 / d$. Thus,
\begin{align*}
\int_{0}^{1} \min \left(\frac{N}{d}, \frac{1}{\| d \alpha \|}\right) d \alpha
 & \ll d \int_{0}^{1 / d} \min \left(\frac{N}{d}, \frac{1}{\| d \alpha \|}\right) d \alpha \\
 & \ll d \left(\int_{0}^{1 / N} \frac{N}{d} \, d \alpha + \int_{1 / N}^{1 / (2 d)} \frac{1}{d \alpha} \, d \alpha \right) \\
 & \ll 1 + \log \left(\frac{N}{2 d}\right)
 \ll \log N.
\end{align*}
Hence,
\[
\int_{0}^{1} \lvert S_{R} (\alpha) \rvert \, d \alpha
 \ll \sum_{d \leq R} \log \left(\frac{R}{d}\right) \log N
 \ll R \log N.
\]
This completes the proof of Lemma \ref{lemma-2}.

The idea here is that $S_{R} (\alpha)$ contributes $N \log R$ in the $L^{2}$ norm but nothing in the $L^{1}$ norm until $R$ is about size $N^{1 / 2}$. Thus, for example,  with $R = N^{1 / 2} / \log N$, we see $S_R(\alpha)$  removes half of the $L^{2}$ norm of $S (\alpha)$, but leaves the $L^{1}$ norm of $S (\alpha)$ unchanged. Hence, the Cauchy--Schwarz inequality gives an improved upper bound for the $L^{1}$ norm of $S (\alpha)$.

\section{Second Proof of Theorem 5} \label{section 7}

We now give another proof of Theorem 5 in which we use part of the circle method to analyze $S (\alpha)$ and directly pull out the major arc contribution. This provides a slightly better approximation for $S (\alpha)$ than $S_{R} (\alpha)$, although the result obtained is asymptotically the same. 

Let $\alpha = a / q + \beta$, where $a / q$ is a fraction with $1 \leq a \leq q$ and $(a, q) = 1$. Suppose $\beta$ is small. Then
\begin{align*}
S \left(\frac{a}{q} + \beta\right)
 & = \sum_{n \leq N} \Lambda (n) e \left(\frac{n a}{q}\right) e (n \beta) \\
 & = \sum_{1 \leq b \leq q} e \left(\frac{b a}{q}\right) \sum_{\substack{n \leq N \\ n \equiv b \tpmod{q}}} \Lambda (n) e (n \beta).
\end{align*}
If $\beta = 0$, by the prime number theorem for arithmetic progressions \cite[Chapter 20, equation (11)]{Davenport2000}, the inner sum is
\[
\sum_{\substack{n \leq N \\ n \equiv b \tpmod{q}}} \Lambda (n)
 = \psi (N; q, b)
 = \operatorname{\boldsymbol{1}}_{(b, q) = 1} \frac{N}{\phi (q)} + O [N \exp (-c_{11} (\log N)^{1 / 2})].
\]
When $\beta$ is small, by partial summation,
\[
\sum_{\substack{n \leq N \\ n \equiv b \tpmod{q}}} \Lambda (n) e (n \beta)
 = \operatorname{\boldsymbol{1}}_{(b, q) = 1} \frac{I (\beta)}{\phi (q)} + \textrm{error}.
\]
Thus, we expect
\[
S \left(\frac{a}{q} + \beta\right)
 = c_{q} (a) \frac{I (\beta)}{\phi (q)} + \textrm{error}
 = \frac{\mu (q)}{\phi (q)} I (\beta) + \textrm{error},
\]
since $c_{q} (a) = \mu (q)$ when $(a, q) = 1$. Hence, we obtain a spike function in each interval which is largest when $q$ is small.  In the circle method we consider the set of all fractions $a / q$ of order $Q$ with $1 \leq q \leq Q$, $1 \leq a \leq q$, and $(a, q) = 1$ and divide the unit interval $[0, 1]$ into nonoverlapping intervals around each fraction $a / q$. Thus, with $\beta = \alpha - a / q$ in \cite{Goldston2000},
\[
J_{Q} (\alpha)
 = \sum_{q \leq Q} \sum_{\substack{1 \leq a \leq q \\ (a, q) = 1}} \frac{\mu (q)}{\phi (q)} I \left(\alpha - \frac{a}{q}\right) \chi_{Q} (\alpha; q, a),
\]
where
\[
\chi_{Q} (\alpha; q, a)
 = \left\{ \begin{array}{ll}
      1 & \mbox{if $\displaystyle \alpha \in \left(\frac{a}{q} - \frac{1}{q (Q + u)}, \frac{a}{q} + \frac{1}{q (Q + v)}\right],\ \  1 \leq u, v \leq q$,} \\
      0 & \mbox{otherwise.}
\end{array}
\right.
\]
We find $J_{Q} (\alpha)$ becomes a better and better approximation for $S (\alpha)$ in $L^{2}$ as $Q$ increases up to size $N^{1 / 2}$ but then becomes a worse approximation as the requirement that the intervals $\chi_{Q} (\alpha; q, a)$ not overlap cuts off the full contribution of each spike function approximation. 

We thus define a new approximation by dropping $\chi_{Q} (\alpha; q, a)$ and letting  each spike $I (\alpha - a / q)$ cut off its contribution just by being small when $\beta$ is not small. 
Thus, we define, replacing $Q$ by $R$,
\[
\mathscr{S}_{R} (\alpha)
 = \sum_{q \leq R} \sum_{\substack{1 \leq a \leq q \\ (a, q) = 1}} \frac{\mu (q)}{\phi (q)} I \left(\alpha - \frac{a}{q}\right).
\]
This is our major arc approximation to $S (\alpha)$, although it does not have arcs. Unlike $J_{Q} (\alpha)$, here we believe that $\mathscr{S}_{R} (\alpha) \to S (\alpha)$ as $R \to N$, but as soon as $R > N^{1 / 2}$ we do not know how to prove this. Thus, we say that $\mathscr{S}_{R} (\alpha)$ is the major arc approximation to $S (\alpha)$ if $R \leq N^{1 / 2}$; but if $R > N^{1 / 2}$, then $\mathscr{S}_{R} (\alpha) - \mathscr{S}_{N^{1 / 2}} (\alpha)$ is the minor arc approximation to $S (\alpha)$. We have
\begin{align*}
\mathscr{S}_{R} (\alpha)
 & = \sum_{q \leq R} \sum_{\substack{1 \leq a \leq q \\ (a, q) = 1}} \frac{\mu (q)}{\phi (q)} \sum_{n \leq N} e (n \alpha) e \left(-\frac{n a}{q}\right) \\
 & = \sum_{n \leq N} \sum_{q \leq R} \frac{\mu (q)}{\phi (q)} \sum_{\substack{1 \leq a \leq q \\ (a, q) = 1}} e \left(-\frac{n a}{q}\right) e (n \alpha) \\
 & = \sum_{n \leq N} \lambda_{R} (n) e (n \alpha), \numberthis \label{eq-19}
\end{align*}
where
\begin{equation} \label{eq-20}
\lambda_{R} (n)
 = \sum_{q \leq R} \frac{\mu (q)}{\phi (q)} c_{q} (-n).
\end{equation}
In comparison to $S (\alpha)$, we expect that $\lambda_{R} (n)$ is an approximation of $\Lambda (n)$. As before, we have
\begin{align*}
\int_{0}^{1} \lvert S (\alpha) - \mathscr{S}_{R} (\alpha) \rvert^{2} \, d \alpha
 & = \sum_{n \leq N} (\Lambda (n) - \lambda_{R} (n))^{2} \\
 & = \sum_{n \leq N} \Lambda (n)^{2} - 2 \sum_{n \leq N} \lambda_{R} (n) \Lambda (n) + \sum_{n \leq N} \lambda_{R} (n)^{2}.
\end{align*}

Let
\[
L (R)
 = \sum_{q \leq R} \frac{\mu (q)^{2}}{\phi (q)}.
\]
By \eqref{eq-9},
\[
L (R)
 = \log R + O (1).
\]
From Lemma \ref{lemma-1} of \cite[equation (2.5)]{Goldston1995} for $1 \leq R \leq N$, and the prime number theorem with remainder we have, for any positive constant $c_{16}$,
\begin{align*}
\sum_{n \leq N} \lambda_{R} (n) \Lambda (n)
 & = \psi (N) L (R) + O (R \log N) \\
 & = N L (R) + O \left( \frac{N}{(\log N)^{c_{16}}}\right) + O (R \log N).
\end{align*}
From Lemma 2 of \cite[equation (2.7)]{Goldston1995}, for $1 \leq R \leq N^{1 / 2}$,
\[
\sum_{n \leq N} \lambda_{R} (n)^{2}
 = N L (R) + O (R^{2}).
\]
Thus, for $R \leq N^{1 / 2}$, by \eqref{eq-9} and \eqref{eq-10},
\begin{align*}
\int_{0}^{1} \lvert S (\alpha) - \mathscr{S}_{R} (\alpha) \rvert^{2} \, d \alpha
 & = \sum_{n \leq N} \Lambda (n)^{2} -NL(R) + O (R^{2})  + O \left( \frac{N}{(\log N)^{c_{16}}}\right) \\
 & = a_{0} (N, R)  + O (R^{2})  + O\left( \frac{N}{(\log N)^{c_{16}}}\right).
\end{align*}
We now replace \eqref{eq-9} by a more precise formula.  Hildebrand \cite[Hilfssatz 2]{Hildebrand1984} proved that, for $x \geq 1$,
\begin{equation} \label{eq-21}
\sum_{\substack{n \leq x \\ (n, k) = 1}} \frac{\mu (n)^{2}}{\phi (n)}
  = \frac{\phi (k)}{k} (\log x + D + g (k)) + O \left(\frac{h (k)}{x^{1 / 2}}\right),
\end{equation}
where
\[
g (k)
 = \sum_{p \mid k} \frac{\log p}{p}, \hspace{1em}
h (k)
 = \sum_{d \mid k} \frac{\mu (d)^{2}}{d^{1 / 2}}
 = \prod_{p \mid k} \left(1 + \frac{1}{p^{1 / 2}}\right).
\]
and
\[
D
 = \gamma + \sum_{p} \frac{\log p}{p (p - 1)} \hspace{1em}
\textrm{($\gamma = $ Euler constant).}
\]
Hence, we obtain a slightly improved estimate over the first proof of 
 \begin{align*}
 \int_{0}^{1} \lvert S (\alpha) - \mathscr{S}_{R} (\alpha) \rvert^{2} \, d \alpha
 & = N \log \left(\frac{N}{R}\right) - N (1 + D) + O (R^{2}) + O \left(\frac{N}{R^{1 / 2}}\right) + O \left(\frac{N}{(\log N)^{A}}\right) \\
 & \leq N\log \left(\frac{N}{R}\right) + O (R^{2}),
\end{align*}
if $N^{1 / 2} / (\log N)^{2} \leq R \leq N^{1 / 2}$.

We complete the proof of Theorem 5 as before by proving the following estimate.

\begin{lemma} \label{lemma-3}
We have
\[
\int_{0}^{1} \lvert \mathscr{S}_{R} (\alpha) \rvert \, d \alpha
 \ll R \log N.
\]
\end{lemma}

To prove Lemma \ref{lemma-3}, substituting the first expression for $c_{q} (-n)$ of \eqref{eq-8} into \eqref{eq-20} yields
\[
\lambda_{R} (n)
 = \sum_{q \leq R} \frac{\mu (q)}{\phi (q)} \sum_{d \mid (q, n)} d \mu \left(\frac{d}{q}\right).
\]
If $d \mid q$, then $q = d m$. We have $\mu (q) = \mu (d m) = \mu (d) \mu (m)$ if $(d, m) = 1$, and $\mu (q) = 0$ if $(d, m) > 1$. Thus, in the nested double sum above, all the nonzero terms when $\mu (q) \neq 0$ are when $(d, m) = 1$, so that $\mu (q / d) = \mu (m) = \mu (q) \mu (d)$. Hence,
\[
\lambda_{R} (n)
 = \sum_{q \leq R} \frac{\mu (q)^{2}}{\phi (q)} \sum_{d \mid (q, n)} d \mu (d)
 = \sum_{\substack{d \leq R \\ d \mid n}} d \mu (d) \sum_{\substack{q \leq R \\ d \mid q}} \frac{\mu (q)^{2}}{\phi (q)}.
\]
Since
\[
\sum_{\substack{q \leq R \\ d \mid q}} \frac{\mu (q)^{2}}{\phi (q)}
 = \sum_{d m \leq R} \frac{\mu (d m)^{2}}{\phi (d m)}
 = \frac{\mu (d)^{2}}{\phi (d)} \sum_{\substack{m \leq R / d \\ (m, d) = 1}} \frac{\mu (m)^{2}}{\phi (m)},
\]
we have
\begin{equation} \label{eq-22}
\lambda_{R} (n)
 = \sum_{\substack{d \leq R \\ d \mid n}} \frac{d \mu (d)}{\phi (d)} \sum_{\substack{m \leq R / d \\ (n, d) = 1}} \frac{\mu (n)^{2}}{\phi (n)}.
\end{equation}
Thus, from \eqref{eq-21},
\begin{align*}
\lambda_{R} (n)
 & = \sum_{\substack{d \leq R \\ d \mid n}} \mu (d) \log \left(\frac{R}{d}\right) + \textrm{smaller terms} \\
 & = \Lambda_{R} (n) + \textrm{smaller terms}.
\end{align*}
This shows that $\Lambda_{R} (n)$ is the first term in the more refined approximation $\lambda_{R} (n)$.

As in the proof for $S_{R} (\alpha)$ in Lemma \ref{lemma-2}, by \eqref{eq-19} and \eqref{eq-22},
\[
\mathscr{S}_{R} (\alpha)
 = \sum_{n \leq N} \sum_{\substack{d \mid n \\ d \leq R}} a (d, R) e (n \alpha)
 = \sum_{d \leq R} a (d, R) \sum_{m \leq N / d} e (m d \alpha),
\]
where
\[
a (d, R)
 = \frac{d \mu (d)}{\phi (d)} \sum_{\substack{m \leq R / d \\ (n, d) = 1}} \frac{\mu (m)^{2}}{\phi (m)}.
\]
Thus, as before,
\begin{align*} 
\int_{0}^{1} \lvert \mathscr{S}_{R} (\alpha) \rvert \, d \alpha
 & \leq \sum_{d \leq R} \lvert a (d, R) \rvert \int_{0}^{1} \min \left(\frac{N}{d}, \frac{1}{\| \alpha d \|}\right) d \alpha \\
 & \ll \sum_{d \leq R} \lvert a (d, R) \rvert \left[1 + \log \left(\frac{N}{2 d}\right)\right] \\
 & \ll \log N \sum_{d \leq R} \log \left(\frac{R}{d}\right)
 \ll R \log N.
\end{align*}
This completes the proof of Lemma \ref{lemma-3}.

\begin{acknowledgments}
\textup{A. Ledoan wishes to express his thanks and appreciation to D.A. Goldston at San Jos\'{e} State University for enlightening discussions and comments on the paper. Ledoan's research is supported by the National Science Foundation under Grant DMS-1852288.}
\end{acknowledgments}

\bibliographystyle{amsplain}

\begin{thebibliography}{00}

\bibitem{Davenport2000} H. Davenport, \textit{Multiplicative number theory}, 3rd ed. (Revised by H.L Montgomery), Graduate Texts in Mathematics 74, Springer-Verlag, New York, 2000.

\bibitem{FunkGoldSengSeng2020} S. Funkhouser, D.A. Goldston, D. Sengupta,  and J. Sengupta, \textit{Prime difference champions}, Discrete mathematics and applications, Springer Optim. Appl., \textbf{165}, 207--236,
Springer, Cham, 2020.

\bibitem{Goldston1984} D.A. Goldston, \textit{The second moment for prime numbers}, Quart. J. Math. Oxford (2) 35 (1984), 153--163.

\bibitem{Goldston1992} D.A. Goldston, \textit{On Bombieri and Davenport's theorem concerning small gaps between primes}, Mathematika, {\bf 39} (1992), 10--17.

\bibitem{Goldston1995} D.A. Goldston, \textit{A lower bound for the second moment of primes in short intervals}, Expo. Math. 13 (1995), 366--376.

\bibitem{Goldston2000} D.A. Goldston, \textit{The major arcs approximation for an exponential sum over primes,} Acta Arith. XCII.2 (2000), 169--179.

\bibitem{GoldstonHuntsNgotiaoco2017} D.A. Goldston, J.Z. Hunts, and T. Ngotiaoco, \textit{The tail of the singular series for the prime pair and Goldbach problems}, Funct. et Approx. Comment. Math. (1) 56 (2017), 117--141.

\bibitem{Graham1978} S. Graham, \textit{An asymptotic estimate related to Selberg's sieve}, J. Number Theory (1) 10 (1978), 83--94.

\bibitem{HardyLittlewood1922} G.H. Hardy and J.E. Littlewood, \textit{Some problems of `Partitio numerorum'; III: On the expression of a number as a sum of primes}, Acta Math. (1) 44 (1922), 1--70. Reprinted as pp. 561--630 in {\it Collected Papers of G.H. Hardy}, Vol. I, Clarendon Press, Oxford University Press, Oxford, 1966.

\bibitem{Heath-Brown1994} D.R. Heath-Brown, \textit{Problem 10,} Tagungsbericht 11/1994, Elementare und Analytische Zahlentheorie, 13.03.1994--19.03.1994, p. 26, Mathematisches Forschungsinstitut Oberwolfach, available at \url{https://www.mfo.de}.

\bibitem{Hildebrand1984} A. Hildebrand, \textit{\"{U}ber die punktweise Konvergenz von Ramanujan-Entwicklungen zahlentheoretischer Funktionen}, Acta Arith. (2) 44 (1984), 109--140.

\bibitem{Ingham1932} A.E. Ingham, \textit{The Distribution of Prime Numbers}, Cambridge Tracts in Mathematics and Mathematical Physics 30, Cambridge University Press, Cambridge, 1932.

\bibitem{Maynard2013} J. Maynard, \textit{Bounded length intervals containing two primes and an almost-prime}, Bull. Lond. Math. Soc. 45 (2013), no. 4, 753--764

\bibitem{Maynard2015} J. Maynard, \textit{Small gaps between primes}, Ann. of Math. (2) 181 (2015), no. 1, 383--413.

\bibitem{MontgomerySoundararajan2002} H.L. Montgomery and K. Soundararajan, \textit{Beyond pair correlation}, Paul Erd\H{o}s and his mathematics, I (Budapest, 1999), 507--514, Bolyai Soc. Math. Stud. 11, Janos Bolyai Math. Soc., Budapest, 2002.

\bibitem{MontgomeryVaughan1973} H.L. Montgomery and R.C. Vaughan, \textit{Error terms in additive prime number theory,} Quart. J. Math. Oxford (2) 24 (1973), 207--216.

\bibitem{MontgomeryVaughan2007} H.L. Montgomery and R.C. Vaughan, \textit{Multiplicative number theory}, Cambridge Studies in Advanced Mathematics 97, Cambridge University Press, Cambridge, 2007.

\bibitem{Polymath2014} D.H.J. Polymath, New equidistribution estimates of Zhang type. Algebra Number Theory 8 (2014), no. 9, 2067--2199.

\bibitem{Polymath2015} D.H.J. Polymath, Variants of the Selberg sieve, and bounded intervals containing many primes. Res. Math. Sci. 1 (2015), Art. 12.

\bibitem{Vaughan1988} R.C. Vaughan, \textit{The $L^{1}$ mean of exponential sums over primes}, Bull. London Math. Soc. 20 (1988), 121--123.

\bibitem{Zhang2014} Y.T. Zhang, \textit{Bounded gaps between primes}, Ann. of Math. (2) 179 (2014), no. 3, 1121--1174.

\end{thebibliography}

\end{document}